\documentstyle[11pt]{article}
\newcommand{\C}{\mbox{C{\llap{{\vrule height1.5ex}\kern.4em}}}\index{C@\mbox{C{\llap{{\vrule height1.3ex}\kern.45em}}}, Notation for the field of complex
numbers}}

\renewcommand{\H}
{{\rm {\mbox{\protect\makebox[.15em][l]{I}H}\index{N@{\protect\makebox[.15em][l]{I}H}, 
Notation for the upper half plane}}}}
\newcommand{\Q}
{{\rm {\mbox{Q{\llap{{\vrule height1.3ex}\kern.5em}}}\index{Q@\mbox{Q{\llap{{\vrule height1.3ex}\kern.5em}}}, 
Notation for the field of rational numbers}}}}
\newcommand{\R}
{{\rm {\mbox{\protect\makebox[.15em][l]{I}R}\index{R@\mbox{\protect\makebox[.15em][l]{I}R}, 
Notation for the field of real numbers}}}}

\newcommand{\ed}[1]{\frac{1}{#1}}

\renewcommand{\Im}{{\rm Im\:}}

\newcommand{\al}{\alpha}

\renewcommand{\phi}{\varphi}

\newcommand{\4}{{\bf{4}}}

\newcommand{\be}{\begin{equation}}
\newcommand{\ee}{\end{equation}}
\newcommand{\bea}{\begin{eqnarray}}
\newcommand{\eea}{\end{eqnarray}}
\newcommand{\beao}{\begin{eqnarray*}}
\newcommand{\eeao}{\end{eqnarray*}}

\newcommand{\lk}{\left(}
\newcommand{\rk}{\right)}

\newcommand{\bi}{\bibitem}

\newcommand{\bT}{\begin{theorem}}
\newcommand{\eT}{\end{theorem}}
\newcommand{\bL}{\begin{lemma}}
\newcommand{\eL}{\end{lemma}}
\newcommand{\bC}{\begin{corollary}}  
\newcommand{\eC}{\end{corollary}}
\newcommand{\bt}{\begin{tabbing} 12345 \= \kill}
\newcommand{\et}{\end{tabbing}}

\newcommand{\abs}{\\[3mm]}

\newtheorem{theorem}{Theorem}

\newtheorem{lemma}{Lemma}
\newtheorem{corollary}{Corollary}

\begin{document}
\begin{center}{
{\LARGE {\bf {A Package on Formal Power Series}}}
\vspace{5mm}\\
{\large {\sc {Wolfram Koepf}}}\\[1mm]
Konrad-Zuse-Zentrum f\"ur Informationstechnik\\
Heilbronner Str.\ 10\\
D-10711 Berlin\\[1mm]
{\tt koepf@zib-berlin.de}
\\[3mm]
Mathematica Journal \4, to appear in May, 1994
}
\end{center}
\vspace{0.5cm}
\begin{center}
{\bf {Abstract:}}
\end{center}%\vspace{0.5cm}
\begin{enumerate}
\item[] {{\small 
Formal Laurent-Puiseux series 
%%of the form
%%$f(x)=\sum\limits_{k=k_0}^{\infty}a_{k}x^{k/n}$
are important in many branches of mathematics. This paper presents a
{\it Mathematica} implementation of algorithms developed by the author
for converting 
between certain classes of functions 
and their equivalent
representing series.
The package {\tt PowerSeries} handles
functions of rational, exponential, and hypergeometric type, and enables
the user to reproduce most of the results of Hansen's extensive table
of series. 
Subalgorithms of independent significance generate
differential equations satisfied by a given function and
recurrence equations satisfied by a given sequence.
\\ }}
\end{enumerate}

\section*{Scope of the Algorithms}
\label{sec:Scope of the algorithm}

A common problem in mathematics is to convert  
an expression involving elementary or special 
functions into its corresponding formal Laurent-Puiseux series of the
form
\be
f(x)=\sum\limits_{k=k_0}^{\infty}a_{k}x^{k/n}
\;.
\label{eq:formalLPS}
\ee
Expressions created from 
algebraic operations on series, such as addition, multiplication,
division,
and substitution, can be handled by finite algorithms if one
truncates the resulting series. These algorithms are implemented in
{\it Mathematica}'s {\tt Series} command.
For example:

{\small
\begin{verbatim}
In[1]:= Series[Sin[x] Exp[x], {x, 0, 5}]

                  3    5
             2   x    x        6
Out[1]= x + x  + -- - -- + O[x]
                 3    30
\end{verbatim}
}\noindent
It is usually much more difficult to find  
the exact formal result, that is, an explicit
formula for the coefficients $a_k$.

This paper presents the package {\tt PowerSeries}, a {\it Mathematica} 
implementation of algorithms developed by the author 
(\cite{Koe92}, \cite{Koe93a}, \cite{Koe93})
for computing Laurent-Puiseux series for three types of functions:
functions of {\sl rational type\/}, which are rational or
have a rational derivative of some order; functions of {\sl exponential
type\/}, which satisfy a homogeneous linear differential equation
with constant coefficients;
and functions of {\sl hypergeometric type\/}, which have a representation
(\ref{eq:formalLPS}) with coefficients satisfying
a recurrence equation of the form
\bea
a_{k+m}&=&R(k)\,a_k\quad\:\mbox{for $k\geq k_0$}
\label{eq:hypergeometric type}
\\
a_{k}&=&A_{k}\quad\quad\quad\;\mbox{for $k=k_0,k_0+1,\ldots,k_0+m-1$}
\nonumber
\eea
for some rational function $R$, 
$A_{k}\in\C$,
$A_{k_0}\neq 0$,
and some positive integer $m$, 
called the {\sl symmetry number} of (the given representation
of) the function.

To find the Laurent-Puiseux expansion of a function $f(x)$ about 
the point $x_0$, one uses the function call 
{\tt PowerSeries[$f(x)$, \{$x$, $x_0$\}]}.
The package also 
extends the built-in
{\tt Series} command 
in the case that its second argument is a list of 
length two:

{\small
\begin{verbatim}
In[2]:= << PowerSeries.m
ps-info: PowerSeries  version 1.02, Mar 07, 1994

In[3]:= Series[Sin[x] Exp[x], {x, 0}]

              k/2  k     k Pi
            2    x  Sin[----]
                         4
Out[3]= Sum[-----------------, {k, 0, Infinity}]
                    k!
\end{verbatim}
}

Our algorithms also handle the converse problem of finding a ``closed form'' 
representation of a given Laurent-Puiseux series (that is, the generating 
function of the sequence $a_k$).  
The standard 
{\it Mathematica} packages {\tt Algebra`SymbolicSum`} 
and {\tt DiscreteMath`RSolve`}
also
deal with this problem, using different approaches. We will discuss
the connection of our package with these packages and other {\it
Mathematica} functions at the end of this introduction.

The most interesting case is formed by the functions of hypergeometric type,
which include 
almost all transcendental elementary 
functions like {\tt x\verb+^+n},
{\tt Exp}, {\tt Log}, 
%%{\tt Sin}, {\tt Sinh}, 
%%{\tt Cos}, {\tt Cosh}, {\tt ArcSin}, {\tt ArcSinh}, 
%%{\tt ArcCos}, {\tt ArcCosh}, {\tt ArcTan}, {\tt ArcTanh},
%%{\tt ArcCot}, {\tt ArcCoth};
and the trigonometric and inverse trigonometric functions;
all kinds of special functions like the Airy functions, 
%%{\tt AiryAi} and 
%%{\tt AiryBi}, 
the Bessel functions,
%% {\tt BesselI}, {\tt BesselJ}, 
%%{\tt BesselY}, and {\tt BesselK}, 
the integral
functions {\tt Erf}, {\tt ExpIntegralEi}, {\tt CosIntegral},
and {\tt SinIntegral}, the associated Legendre functions
{\tt LegendreP} and
{\tt LegendreQ}, the hypergeometric functions
{\tt HypergeometricU},
{\tt Hypergeometric0F1},
{\tt Hypergeometric1F1}, as well as %\linebreak
{\tt Hypergeometric2F1}, 
the orthogonal polynomials
{\tt JacobiP}, \linebreak
{\tt GegenbauerC}, {\tt ChebyshevT},
{\tt ChebyshevU}, {\tt LegendreP},
{\tt LaguerreL}, and \linebreak
{\tt HermiteH}; and many more
functions.
%Some of the given examples have so-called logarithmic
%singularities which, however, can be covered by the given approach.

Our algorithms also handle the following special functions.
Since they are not built in to 
{\it Mathematica}, 
we have included the necessary
definitions in our package: the Bateman functions {\tt Bateman[n,x]}
(\cite{AS}, (13.6), and \cite{KoeBateman}); the Hankel functions
{\tt Hankel1[n,x]} and {\tt Hankel2[n,x]} (\cite{AS}, (9.1));
the Kummer and Whittaker functions {\tt KummerM[a,b,x]}, {\tt KummerU[a,b,x]},
{\tt WhittakerM[a,b,x]}, {\tt WhittakerW[a,b,x]} (\cite{AS}, (13.1)); 
the Struve functions
{\tt StruveH[n,x]} and {\tt StruveL[n,x]} (\cite{AS}, Chapter 12);
the hypergeometric functions
{\tt Hypergeometric1F0[a,x]}, {\tt Hypergeometric2F0[a,b,x]};
the iterated integrals of the complementary error function {\tt Erfc[n,x]}
(\cite{AS}, (7.2)); the Abramowitz functions {\tt Abramowitz[n,x]}
(\cite{AS}, (27.5)); and the repeated integrals of the normal probability 
integral {\tt NormalIntegral[n,x]} (\cite{AS}, (26.2.41)).
Informations about these functions 
can be obtained by 
using {\it Mathematica}'s {\tt ?} feature
after the package is loaded.

Our algorithms
depend strongly on the fact that any
function of hypergeometric type
satisfies a {\it simple differential equation}, 
i.\ e., a homogeneous
linear differential equation with polynomial coefficients
(see \cite{Koe92}, Theorems 2.1, and 8.1). 
One of our
main algorithms finds this differential equation for (practically)
any expression for which such an equation  exists.

On the other hand, it is similarly important that almost every
function that we may write down satisfies a simple differential equation.
All functions that can be
constructed algebraically from the functions mentioned above by
addition, multiplication, and 
composition with rational functions and rational powers, satisfy such a
differential equation
(\cite{Sta}, \cite{Koediffgl}). 

The conversion of functions 
to and from their representing
Laurent-Puiseux expansions 
is described in detail in
(\cite{Koe92}, \cite{Koe93}).
Here, we give a brief description
and some examples for each of the types of functions considered.

To find the Laurent-Puiseux expansion of an expression $f(x)$,
the function {\tt PowerSeries} performs the following steps.
First, a 
homogeneous linear differential equation with polynomial
coefficients for $f$ is generated. (If no low-order
differential equation exists, this may take some time.)
This differential equation is converted to an equivalent recurrence
equation for $a_k$, an easy and fast procedure.
If the recurrence equation is of the hypergeometric type, the 
coefficients can be calculated explicitly from a finite number of
initial values. 
(Finding the initial values may require the calculation of limits. 
This procedure may
be slow, or may fail.)

As an example of a function of rational type, 
let's take $f(x)=\arctan x $. Here, $f'(x)=\ed{1+x^2}$ is rational,
so we apply a special algorithm for the rational case, and integrate. 
First, the complex partial fraction decomposition is derived:
\[
f'(x)=\frac{1}{1+x^2}=\ed 2\lk\ed{1+ix}+\ed{1-ix}\rk \;,
\]
from which the coefficients
$b_k$ of the derivative $F'(x)=\sum\limits_{k=0}^\infty b_k x^k$
can be deduced:
\[
b_k=\ed 2\lk i^k+(-i)^k\rk=\frac{i^k}{2}\lk 1+(-1)^k\rk
\;,
\]
Finally, after an integration, we get
\[
f(x)=
\sum _{k=0}^{\infty }{\frac {i^k\,{\left ( 1+(-1)^{k}\right )}}
{2\,\left (k+1\right )}}\,x^{k+1}
\;.
\]
To illustrate the 
hypergeometric case,
we again take $f(x)=\arctan x $. We have
$f'(x)=\ed{1+x^2}$ and $f''(x)=-{\frac {2\,x}{\left
(1+x^{2}\right )^{2}}}$, so we arrive at the differential equation
\[
(1+x^{2})\,f''(x)+2\,x \,f'(x)=0 \;.
\]
This differential equation can be converted to the recurrence equation
\be
\left (k+1\right )\left (k+2\right )\,a_{k+2}+k\left (k+1\right )\,a_k=0
\label{eq:arctanRE}
\ee
by the rule-based function {\tt detore[DE,F,x,a,k]}, which 
takes as input the left hand side of the differential equation for the function 
{\tt F[x]} and computes the left hand side of
the recurrence equation 
for the coefficients {\tt a[k]}:

{\small
\begin{verbatim}
detore[g_ + h_, F_, x_, a_, k_]:= 
        detore[g, F, x, a, k] + detore[h, F, x, a, k]
detore[c_ * g_, F_, x_, a_, k_]:= c * detore[g, F, x, a, k] /; 
        FreeQ[c, x] && FreeQ[c, F]
detore[Derivative[k0_][F_][x_], F_, x_, a_, k_]:=
        Pochhammer[k+1, k0] * a[k+k0]
detore[F_[x_], F_, x_, a_, k_]:= a[k]
detore[x_^j_. * Derivative[k0_][F_][x_], F_, x_, a_, k_]:=
        Pochhammer[k+1-j, k0] * a[k+k0-j]
detore[x_^j_. * F_[x_], F_, x_, a_, k_]:= a[k-j]
\end{verbatim}
}\noindent

{\small
\begin{verbatim}
In[4]:= Simplify[detore[2*x*f'[x] + f''[x] + x^2*f''[x], f, x, a, k]]

Out[4]= (1 + k) (k a[k] + 2 a[2 + k] + k a[2 + k])
\end{verbatim}
}\noindent
As a last step, this recurrence equation is solved, resulting in the
representation
\[
f(x)=\sum _{k=0}^{\infty }{\frac {\left (-1\right)^{k}}{2\,k+1}}\,x^{2\,k+1}
\;.
\]
Note that in most cases the hypergeometric-type result, if applicable,
is the simplest one.

Finally, let us take $f(x)=\sin x\,e^x$, the function  
whose series was calculated above.
Then $f'(x)=e^{x}\Big(\sin x
+\cos x \Big)$ and $f''(x)=2e^{x}\cos x $. 
If we look for an equation 
of the form
\[
f''+A_{1}f'+A_{0}f
=2e^{x}\cos\:x+A_{1}e^{x}\Big(\!\sin\:x+\cos\:x\!\Big)+A_{0}e^{x}\sin\:x
=0
\;,
\]
we may expand and recombine rationally dependent expressions to get
\[
f''+A_{1}f'+A_{0}f
=\Big( 2+A_{1} \Big) e^{x}\cos\:x +\Big(A_{1} +A_{0}\Big)  e^{x}\sin\:x
=0
\;.
\]
This equation must be valid for all $x$, so 
the coefficients of $e^{x}\cos\:x$ and
$e^{x}\sin\:x$ must vanish. Hence $A_{1}=-2$, $A_{0}=2$, and 
we have the differential equation for $f$:
\[
f''-2f'+2f=0\;.
\]
This differential equation has constant coefficients, so $f$ is of exponential 
type. For $b_n= n! \,a_n$, we have the recurrence equation
\[
b_{n+2} -2b_{n+1}+2b_{n}=0
\]
with initial values
\[
b_0=a_0=f(0)=0\;\quad\quad{\rm and}\quad\quad
b_1=a_1=\frac{f'(0)}{1!}=1
\;,
\]
which generates the solution
\beao
a_n
&=&
\frac{b_n}{n!} = \frac{1}{n!}\frac{(1+i)^n-(1-i)^n}{2i}
\\& = &
\ed{n!}\Im (1+i)^n = \ed{n!}\Im \lk\sqrt 2 e^{i\frac{\pi}{4}}\rk ^n
\\& = &
\ed{n!} 2^{\frac{n}{2}} \sin {\frac{n \pi}{4}}
\;,
\eeao
Hence
\[
f(x)=\sum\limits_{n=0}^\infty \ed{n!} 2^{\frac{n}{2}}
\sin {\frac{n \pi}{4}} \,x^n\;.
\]
We will present more examples to demonstrate our {\it Mathematica} implementation 
later in the paper. Additional examples are also given in \cite{Koe93a}.

To find the generating function of a given sequence, 
the process can be reversed. The function 
{\tt Convert} 
provides the following procedure to calculate the closed-form
representation of $\sum\limits_{k=k_0}^\infty a_k x^{k/n}$.
First, 
a homogeneous linear recurrence equation with polynomial
coefficients for $a_k$ is generated. (If no low order
recurrence equation exists, this may take some time.)
This recurrence equation is converted to an equivalent differential
equation for the generating function, also an easy and fast
procedure.
Finally, the differential equation is solved using appropriate initial values.
(This step depends heavily on {\it Mathematica}'s built-in differential
equation solver {\tt DSolve}. It often takes a lot of time to solve
the differential equation, or may also fail. Moreover,
the calculation of limits is needed to
find the initial values. This part of the
procedure may be slow or fail, as well.)

As an example,
we try to find the generating function
\[
\sum_{k=0}^\infty \frac{(-1)^k}{2k+1}\,x^k
\]
of the sequence $a_k=\frac{(-1)^k}{2k+1}$. Taking a shift, we get the
relation $a_{k+1}=\frac{(-1)^{k+1}}{2k+3}$ so that the recurrence
equation
\[
(2k+3) a_{k+1}+(2k+1) a_{k}=0
\]
is established for $k\geq 0$. If we multiply this 
equation by the factor $(k+1)$,
the resulting recurrence equation
\be
(k+1) (2k+3) a_{k+1}+(k+1)(2k+1) a_{k}=0
\label{eq:LHS}
\ee
holds for all integers $k$. The back-substitution can be done by the
{\it Mathematica} procedure {\tt retode[re,n,f,x]}, defined by these rules:

{\small
\begin{verbatim}
theta[f_, x_]:= x D[f, x]
retode[eq1_ + eq2_, k_, f_, x_]:= 
        retode[eq1, k, f, x] + retode[eq2, k, f, x]
retode[c_ * eq_, k_, f_, x_]:= c * retode[eq, k, f, x] /;
        (FreeQ[c, k] && FreeQ[c, f] && FreeQ[c, x])
retode[a[k_ + m_.], k_, f_, x_]:= f[x]/x^m
retode[k^j_. * a[k_ + m_.], k_, f_, x_]:=
        theta[retode[k^(j-1) * a[k+m], k, f, x], x]
retode[p_ * a[k_ + m_.], k_, f_, x_]:= 
        retode[Expand[p * a[k+m]], k, f, x] /; PolynomialQ[p]
\end{verbatim}
}
\noindent
Applied to the left hand side of (\ref{eq:LHS}), {\tt retode} gives:

{\small
\begin{verbatim}
In[5]:= Simplify[retode[(k+1)(2k+3)a[k+1] + (k+1)(2k+1)a[k], k, f, x]]

Out[5]= f[x] + 3 f'[x] + 5 x f'[x] + 2 x f''[x] +

        2
>    2 x  f''[x]
\end{verbatim}
}
\noindent
The initial value problem
\[
(2x+2x^2) f''+(3+5x)f'+f=0\;,
\]
\[
f(0)=a_0=1\;,\quad f'(0)=a_1=-\ed 3\;,
\]
has the solution
\[
f(x)=\frac{\arctan\sqrt x}{\sqrt x}\;.
\]
Note that {\it Mathematica} Version 2.0 gives the correct result:

{\small
\begin{verbatim}
In[6]:= Convert[Sum[(-1)^k/(2k+1)*x^k, {k, 0, Infinity}], x]

         ArcTan[Sqrt[x]]
Out[6]= ---------------
             Sqrt[x]
\end{verbatim}
}
\noindent
whereas the {\tt DSolve} command of {\it Mathematica} Version 2.2
gives a wrong result.

Now we want to give some remarks on the connection of our package with
several standard {\it Mathematica} packages.

In several instances, {\tt PowerSeries} depends on the calculation of
limits, so it may be helpful to load the package {\tt Calculus`Limit`}.
Since the {\tt Convert} function of {\tt PowerSeries} depends heavily on
the {\tt DSolve} command, it may also be helpful to load 
{\tt Calculus`DSolve`} (Version 2.2). 
In this case, the incorrect result mentioned above does not occur.

The aim of the package {\tt Algebra`SymbolicSum`}  
is to give closed
form representations of sums, in particular, sums of the generating
function type.
Thus, it is not surprising that many of the conversions that can be obtained
by our {\tt Convert} function are also obtained by the {\tt SymbolicSum} 
package. 
Among the expressions that cannot be handled by 
the {\tt SymbolicSum} package 
are the sums {\tt conv[4]} and
{\tt conv[5]} defined below, and, in general, 
sums involving special functions. 
Typically, our implementation is slower as it always uses {\tt DSolve},
which takes a lot of time. The advantage is that the
procedure forms an algorithmic treatment of the problem, and is not
based on heuristics.

%If you have both, {\tt SymbolicSum} and {\tt PowerSeries} loaded,
%and issue a {\tt Convert} command, then
%
%This may be, probably, a good combination,
%although sometimes {\tt SymbolicSum}'s output contains more
%see e.\ g.\ our example {\tt conv[3]} in
%\S~\ref{sec:Examples on the calculation of generating functions}.

The purpose of the package {\tt DiscreteMath`RSolve`}
is to solve recurrence equations.
It uses generating function techniques and defines the function
{\tt GeneratingFunction} to find the generating
function of a recurrence equation. The procedure uses
heuristic techniques, and in particular treats hypergeometric
recurrence equations (equation (\ref{eq:hypergeometric type}) with $m=1$), but
the general case ($m>1$) is not implemented. Similar restrictions apply
to the {\tt RSolve} command.

We give an example:
Our function 
{\tt SimpleRE} 
generates the
recurrence equation (\ref{eq:arctanRE})
for the coefficients of the series representation
of a function:

{\small
\begin{verbatim}
In[7]:= re=SimpleRE[ArcTan[x], x]

Out[7]= k (1 + k) a[k] + (1 + k) (2 + k) a[2 + k] == 0
\end{verbatim}
}
\noindent
We load the package  and try to solve the recurrence equation:

{\small
\begin{verbatim}
In[8]:= << DiscreteMath`RSolve`

In[9]:= RSolve[{re, a[0]==0, a[1]==1}, a[k], k]

DSolve::dnim:
   Built-in procedures cannot solve this differential
     equation.

Out[9]= RSolve[{re, a[0] == 0, a[1] == 1}, a[k], k]
\end{verbatim}
}
\noindent
(This is the {\it Mathematica} Version 2.2 result. In Version 2.0,
another error message occurs.)
We see that {\tt RSolve} tries to solve the differential equation
for the generating function of the sequence $a_k$ given by the recurrence 
equation 
\be
k \,(1 + k)\, a_k + (1 + k)\, (2 + k) \,a_{k+2} = 0
\; .
\label{eq:simpler one}
\ee
This problem is too hard. In our treatment, the recurrence
equation is replaced by the simpler equation
\[
k \,(1 + k)\, a_k + (1 + k)\, (2 + k) \,a_{k+1} = 0
\]
whose generating function is easier to find. Further, our
function {\tt PowerSeries} solves 
the recurrence equation (\ref{eq:simpler one}) by a direct application
of the hypergeometric theory, and splits the series into two parts
(odd and even), one of which disappears as a result of
a vanishing initial value:

{\small
\begin{verbatim}
In[10]:= Series[ArcTan[x], {x, 0}]

                 k  1 + 2 k
             (-1)  x
Out[10]= Sum[--------------, {k, 0, Infinity}]
                1 + 2 k
\end{verbatim}
}
\noindent

\section*{Generation of Differential Equations}
\label{sec:Examples of the generation of differential equations}

In this section, we present some 
examples of the
function {\tt SimpleDE}, 
which generates the 
homogeneous linear differential equation of least order satisfied by 
a given function.

{\small
\begin{verbatim}
In[11]:= de[1] = SimpleDE[x^n * E^(alpha x), x]

Out[11]= (-n - alpha x) F[x] + x F'[x] == 0

In[12]:= de[2] = SimpleDE[((1+x)/(1-x))^n, x]

Out[12]= 2 n F[x] + (-1 + x) (1 + x) F'[x] == 0

In[13]:= de[3] = SimpleDE[ArcSin[x^5], x]

               10                 11
Out[13]= (4 + x  ) F'[x] + (-x + x  ) F''[x] == 0

In[14]:= de[4] = SimpleDE[ArcSin[x], x]

                          2
Out[14]= x F'[x] + (-1 + x ) F''[x] == 0
\end{verbatim}
}\noindent
As noted above, 
we can create new functions 
satisfying this kind of differential equation 
by addition, multiplication, and by the
composition with rational functions and rational powers:

{\small
\begin{verbatim}
In[15]:= de[5] = SimpleDE[ArcSin[x]^3, x]

                            2
Out[15]= x F'[x] + (-4 + 7 x ) F''[x] +

                 2   (3)             2 2  (4)
>     6 x (-1 + x ) F   [x] + (-1 + x )  F   [x] == 0

In[16]:= de[6] = SimpleDE[Sin[x]^5, x, 6]

                                     (4)
Out[16]= 225 F[x] + 259 F''[x] + 35 F   [x] +

       (6)
>     F   [x] == 0

In[17]:= de[7] = SimpleDE[Exp[alpha x] * Sin[beta x], x]

               2       2
Out[17]= (alpha  + beta ) F[x] - 2 alpha F'[x] +

>     F''[x] == 0
\end{verbatim}
}\noindent
By default, {\tt SimpleDE} searches for a simple differential equation of order
up to five. This setting can be changed by entering the order as
a third argument. In case one knows
a priori that a differential equation of the given type exists, one may choose 
the order {\tt Infinity}.
Note that there is no simple differential equation of order less than five
for $\sin^5 x$, so {\tt SimpleDE[Sin[x]\verb+^+5,x]} fails.

Our implementation handles
orthogonal polynomials and special functions 
using the general method described 
in \cite{Koeortho}. 
We can derive the usual
differential equations:

{\small
\begin{verbatim}
In[18]:= de[8] = SimpleDE[LaguerreL[n, x], x]

Out[18]= n F[x] + (1 - x) F'[x] + x F''[x] == 0

In[19]:= de[9] = SimpleDE[ChebyshevT[n, x], x]

            2                          2
Out[19]= -(n  F[x]) + x F'[x] + (-1 + x ) F''[x] == 0

In[20]:= de[10] = SimpleDE[BesselY[n, x], x]

            2    2                    2
Out[20]= (-n  + x ) F[x] + x F'[x] + x  F''[x] == 0

In[21]:= de[11] = SimpleDE[AiryAi[x], x]

Out[21]= -(x F[x]) + F''[x] == 0
\end{verbatim}
}\noindent
By the statement above, we know that much more complicated functions 
satisfy differential equations:

{\small
\begin{verbatim}
In[22]:= de[12] = SimpleDE[Exp[alpha x] * BesselI[n, x], x]

            2              2        2  2
Out[22]= (-n  - alpha x - x  + alpha  x ) F[x] +

                    2           2
>     (x - 2 alpha x ) F'[x] + x  F''[x] == 0

In[23]:= de[13] = SimpleDE[Sin[m x] * BesselJ[n, x], x]

           2      4      6    2      2  2       2  2
Out[23]= (n  - 5 n  + 4 n  - x  - 2 m  x  + 22 n  x  - 
 
             2  2  2       4  2       2  4  2    4
>        10 m  n  x  - 12 n  x  + 12 m  n  x  - x  + 
\end{verbatim}

%%\pagebreak
}\noindent
{\small
\begin{verbatim}
            2  4      4  4       2  4       2  2  4
>        6 m  x  - 5 m  x  + 12 n  x  - 24 m  n  x  + 
 
             4  2  4      6       2  6       4  6
>        12 m  n  x  - 4 x  + 12 m  x  - 12 m  x  + 
 
            6  6
>        4 m  x ) F[x] + 
 
           3      2  3      2  3       2  2  3
>     (-8 x  - 4 m  x  + 8 n  x  + 40 m  n  x  - 
 
            5      4  5
>        8 x  + 8 m  x ) F'[x] + 

           2       2  2      4  2      4      2  4
>     (-2 x  + 10 n  x  - 8 n  x  + 2 x  - 6 m  x  + 
 
             2  4      6      4  6
>        16 n  x  - 8 x  + 8 m  x ) F''[x] + 
 
           3       2  3      5      2  5   (3)
>     (-4 x  + 16 n  x  - 8 x  + 8 m  x ) F   [x] + 
 
       4          2      2      2  2   (4)
>     x  (-1 + 4 n  - 4 x  + 4 m  x ) F   [x] == 0

In[24]:= de[14] = SimpleDE[LegendreP[n, x]^2, x]

                      2
Out[24]= (-4 n x - 4 n  x) F[x] +

                     2      2        2      2  2
>     (-2 + 4 n + 4 n  + 6 x  - 4 n x  - 4 n  x )

                          2
>      F'[x] + 6 x (-1 + x ) F''[x] +

             2 2  (3)
>     (-1 + x )  F   [x] == 0
\end{verbatim}
}\noindent
(The calculation of these last two results
is very time consuming!)

The functions
\[
F_n^{(\al)}(x):=e^{-x}\,L_n^{(\al)}(2x)
\;,
\]
where $L_n^{(\al)}$ denotes the generalized Laguerre polynomial,
were studied by the author in \cite{KoeBateman}.
For these functions, we immediately deduce the differential equation:

{\small
\begin{verbatim}
In[25]:= de[15] = SimpleDE[E^(-x) * LaguerreL[n, alpha, 2x], x]

Out[25]= (1 + alpha + 2 n - x) F[x] +

>     (1 + alpha) F'[x] + x F''[x] == 0
\end{verbatim}
}\noindent
We may observe that 
the coefficient of
$F'$ in the differential equation vanishes 
only in the case $\al=-1$, which  
is studied in \cite{KoeBateman}.

\section*{Conversion of Simple Differential Equations to
Recurrence Equations}
\label{sec:Examples on the conversion of simple differential equations
to recurrence equations}

The function {\tt DEtoRE} converts a simple differential
equation into a recurrence 
equation for the coefficients of a corresponding Laurent-Puiseux 
expansion. The output 
is given in terms of the symbols $a$ and $k$ representing the
coefficient sequence $a_k$.  
As examples, we convert some of the differential equations that
were derived above.

The example

{\small
\begin{verbatim}
In[26]:= re[1] = DEtoRE[de[2], F, x]

Out[26]= k a[k] + 2 n a[1 + k] - (2 + k) a[2 + k] == 0
\end{verbatim}
}\noindent
shows that the function $\lk\frac{1+x}{1-x}\rk^n$ is not of hypergeometric
type, whereas the example

{\small
\begin{verbatim}
In[27]:= re[2] = DEtoRE[de[3], F, x]

          2
Out[27]= k  a[k] - (5 + k) (10 + k) a[10 + k] == 0
\end{verbatim}
}\noindent
shows that the function $\arcsin\:(x^5)$ is of hypergeometric type with
symmetry number $m=10$. Similarly,
the Airy function is of hypergeometric type with symmetry number $m=3$:

{\small
\begin{verbatim}
In[28]:= re[3] = DEtoRE[de[11], F, x]

Out[28]= -a[k] + (2 + k) (3 + k) a[3 + k] == 0
\end{verbatim}
}\noindent
Among the functions of the form 
$e^{\al x}\, I_n (x)$,
where $\al\in\R$ and $I$ is the modified Bessel function of the first kind,
only the cases $\al=0$ and $\al=\pm 1$ represent functions of 
hypergeometric type, as is seen from the recurrence equation

{\small
\begin{verbatim}
In[29]:= re[4] = DEtoRE[de[12], F, x]

Out[29]= (-1 + alpha) (1 + alpha) a[k] -

>     alpha (3 + 2 k) a[1 + k] +

>     (2 + k - n) (2 + k + n) a[2 + k] == 0
\end{verbatim}
}\noindent
Similarly, the
functions of the form 
$\sin\:(m x)\, J_m (x)$, 
where $(m\in\R)$ and $J$ is the Bessel function of the first kind,
are of hypergeometric type
only for the cases $m=0$ and $m=\pm 1$, 
as is seen from the recurrence equation

{\small
\begin{verbatim}
In[30]:= re[5] = DEtoRE[de[13], F, x]

                   3        3
Out[30]= 4 (-1 + m)  (1 + m)  a[k] +

>     (-1 + m) (1 + m)

                       2       2         2      2  2
>      (33 + 32 k + 8 k  + 27 m  + 32 k m  + 8 k  m  -

             2       2  2
>        12 n  + 12 m  n ) a[2 + k] +

                           2       3      4        2
>     (-297 - 402 k - 210 k  - 48 k  - 4 k  + 198 m  +

                2        2  2       3  2      4  2
>        362 k m  + 206 k  m  + 48 k  m  + 4 k  m  +

              2          2       2  2        2  2
>        246 n  + 120 k n  + 16 k  n  + 150 m  n  +

               2  2       4       2  4
>        40 k m  n  - 12 n  + 12 m  n ) a[4 + k] +

>     (5 + k - n) (6 + k - n) (5 + k + n) (6 + k + n)

>      (-1 + 2 n) (1 + 2 n) a[6 + k] == 0
\end{verbatim}
}\noindent
The function {\tt SimpleRE} combines the two steps {\tt SimpleDE} and
{\tt DEtoRE} to calculate
the recurrence equation of a Laurent-Puiseux series representation of 
a function.

\section*{Calculation of Laurent-Puiseux Series}
\label{sec:Examples on the calculation of Laurent-Puiseux series}

In this section, we present some of the results of the procedure 
{\tt PowerSeries}. Obviously, the program can handle the standard power series 
representations of the elementary functions:

{\small
\begin{verbatim}
In[31]:= ps[1] = PowerSeries[E^x, x]

              k
             x
Out[31]= Sum[--, {k, 0, Infinity}]
             k!
\end{verbatim}

%%\pagebreak
}\noindent
{\small
\begin{verbatim}
In[32]:= ps[2] = PowerSeries[Sin[x], x]

                 k  1 + 2 k
             (-1)  x
Out[32]= Sum[--------------, {k, 0, Infinity}]
               (1 + 2 k)!


In[33]:= ps[3] = PowerSeries[ArcSin[x], x]

              1 k  1 + 2 k       2
             (-)  x        (2 k)!
              4
Out[33]= Sum[---------------------, {k, 0, Infinity}]
                  2
                k!  (1 + 2 k)!

\end{verbatim}
}\noindent
Here is a series expanded about the point $x=1$:
{\small
\begin{verbatim}
In[34]:= ps[4] = PowerSeries[Log[x], {x, 1}]

                 k         1 + k
             (-1)  (-1 + x)
Out[34]= Sum[-------------------, {k, 0, Infinity}]
                    1 + k
\end{verbatim}

It may come as a surprise that the following functions are of
hypergeometric type:

{\small
\begin{verbatim}
In[35]:= ps[5] = PowerSeries[Exp[ArcSin[x]], x]

              k  2 k         5            2
             4  x    Product[- - 2 jj + jj , {jj, k}]
                             4
Out[35]= Sum[----------------------------------------,
                              (2 k)!

>     {k, 0, Infinity}] +

          k  1 + 2 k         1          2
         4  x        Product[- - jj + jj , {jj, k}]
                             2
>    Sum[------------------------------------------,
                         (1 + 2 k)!

>     {k, 0, Infinity}]

In[36]:= ps[6] = PowerSeries[Exp[ArcSinh[x]], x]

                    1  k  2 k
                 (-(-))  x    (2 k)!
                    4
Out[36]= x + Sum[-------------------, {k, 0, Infinity}]
                                2
                    (1 - 2 k) k!
\end{verbatim}
}\noindent
If one is interested in the intermediate calculations, one may give the input 
{\tt psprint}, which sets the global variable {\tt PSPrintMessages} to
{\tt True}:

{\small
\begin{verbatim}
In[37]:= psprint

In[38]:= ps[7] = PowerSeries[E^x - 2 E^(-x/2) Cos[Sqrt[3]x/2 - Pi/3], x]
ps-info: PowerSeries  version 1.02, Mar 07, 1994
ps-info: 3 step(s) for DE:
                  (3)
         -F[x] + F   [x] == 0
ps-info: RE for all k >= 0:
         a[3 + k] = a[k]/((1 + k)*(2 + k)*(3 + k))
ps-info: function of hypergeometric type
ps-info: a[0] = 0
ps-info: a[1] = 0
                3
ps-info: a[2] = -
                2
ps-info: a[3] = 0
ps-info: a[4] = 0

                        2 + 3 k
             9 (1 + k) x
Out[38]= Sum[------------------, {k, 0, Infinity}]
                 (3 + 3 k)!
\end{verbatim}
}\noindent
Note that the differential equation in this example is extremely simple.

We now suppress intermediate results using the command 
{\tt nopsprint}, which sets {\tt PSPrintMessages} to {\tt False}.

{\small
\begin{verbatim}
In[39]:= nopsprint
\end{verbatim}
}\noindent

Here is a well-known closed formula for 
the generating function $f(x)=\sum\limits_{k=0}^\infty a_k x^k$
of the Fibonacci numbers $a_n$ defined by the recurrence
\[
a_{n+1}=a_n+a_{n-1}\;,
\quad\quad\quad
a_0=0\;,\quad a_1=1\;.
\]

{\small
\begin{verbatim}
In[40]:= ps[8] = PowerSeries[x/(1 - x - x^2), x]

               k       1       k        -2      k   k
             (2  (------------)  - (-----------) ) x
                  -1 + Sqrt[5]      1 + Sqrt[5]
Out[40]= Sum[----------------------------------------, 
                             Sqrt[5]
 
>    {k, 0, Infinity}]
\end{verbatim}
}\noindent
Moreover, 
we can use 
{\tt SimpleRE} to generate the Fibonacci recursion:

{\small
\begin{verbatim}
In[41]:= re[6] = SimpleRE[x/(1 - x - x^2), x, a, k]

Out[41]= (1 + k) a[k] + (1 + k) a[1 + k] -

>     (1 + k) a[2 + k] == 0
\end{verbatim}
}\noindent
Note that the common factor $(1+k)$ ensures that the
recurrence equation holds for all integers $k$.

Our implementation covers functions that correspond to hypergeometric type 
Laurent-Puiseux series, rather than just power series. For examples, see
\cite{Koe93a}.

Here, we want to return to the case of special functions. 
In the previous section, 
we discovered that the functions
$e^{x}\, I_n (x)$ and $\sin\:(x)\, J_n (x)$
are of hypergeometric type. For $n=0$ and $1$, we get the following series
representations:

{\small
\begin{verbatim}
In[42]:= ps[9] = PowerSeries[E^x * BesselI[0, x], x]

              1 k  k
             (-)  x  (2 k)!
              2
Out[42]= Sum[--------------, {k, 0, Infinity}]
                    3
                  k!

In[43]:= ps[10] = PowerSeries[E^x * BesselI[1, x], x]

              1 k  1 + k
             (-)  x      (1 + 2 k)!
              2
Out[43]= Sum[----------------------, {k, 0, Infinity}]
                    2
                  k!  (2 + k)!
\end{verbatim}

%%\pagebreak
}\noindent
{\small
\begin{verbatim}
In[44]:= ps[11] = PowerSeries[Sin[x] * BesselJ[0, x], x]

                1  k  1 + 2 k
             (-(-))  x        (1 + 4 k)!
                4
Out[44]= Sum[---------------------------,
                                  2
                 (2 k)! (1 + 2 k)!

>    {k, 0, Infinity}]

In[45]:= ps[12] = PowerSeries[Sin[x] * BesselJ[1, x], x]

                1  k  2 + 2 k
             (-(-))  x        (3 + 4 k)!
                4
Out[45]= Sum[---------------------------,
                          2
              2 (1 + 2 k)!  (3 + 2 k)!

>    {k, 0, Infinity}]
\end{verbatim}
}\noindent
% COMMENT In[37] braucht ca. 5 Minuten
% COMMENT In[38] braucht ca. 5 Minuten
Many more functions constructed from 
the given special functions 
are of hypergeometric type  (see, for example, Hansen's
extensive table of series \cite{Han}),
and their Laurent-Puiseux expansions can therefore be found by
{\tt PowerSeries}.

\section*{Calculation of Generating Functions}
\label{sec:Examples on the calculation of generating functions}

The function {\tt Convert} calculates
generating functions of series:

{\small
\begin{verbatim}
In[46]:= conv[1] = Convert[Sum[(2k)!/k!^2 x^k, {k, 0, Infinity}], x]

               1
Out[46]= -------------
         Sqrt[1 - 4 x]

In[47]:= conv[2] = Convert[Sum[k!^2/(2k)! x^k, {k, 0, Infinity}], x]

                                  Sqrt[x]
                 4 Sqrt[x] ArcSin[-------]
           4                         2
Out[47]= ----- + -------------------------
         4 - x                 3/2
                        (4 - x)
\end{verbatim}
}\noindent
Again, with {\tt psprint} we get information about the intermediate steps:

{\small
\begin{verbatim}
In[48]:= psprint

In[49]:= conv[3] = Convert[Sum[k!/(2k)! x^k, {k, 0, Infinity}], x]

ps-info: PowerSeries  version 1.02, Mar 07, 1994
ps-info: 1 step(s) for RE.
         a[-1 + k] + 2*(1 - 2*k)*a[k] == 0
ps-info: DE: 
         F[x] + (-2 + x) F'[x] + (-4*x) F''[x] == 0
ps-info: Trying to solve DE ...
ps-info: DSolve`t encountered
ps-info: DSolve computes 
         C[3] + E^(x/4)*x^(1/2)*
 
>     (C[1] + (C[2]*(-2 - 
 
>            E^(x/4)*Pi^(1/2)*x^(1/2)*Erf[x^(1/2)/2]))/
 
>        (E^(x/4)*x^(1/2)))
ps-info: expression rearranged:
         E^(x/4)*x^(1/2)*C[1] + C[3] + 
 
>    E^(x/4)*x^(1/2)*C[2]*Erf[x^(1/2)/2]
ps-info: Calculation of initial values...
ps-info: C[3] == 1
ps-info: C[1] == 0
ps-info: C[2]/Pi^(1/2) == 1/2

              x/4                      Sqrt[x]
             E    Sqrt[Pi] Sqrt[x] Erf[-------]
                                          2
Out[49]= 1 + ----------------------------------
                             2
\end{verbatim}
}\noindent

Here are two generating function examples involving special functions.
(Note that {\tt Convert} 
can only be successful if the
solution is representable by an elementary function; otherwise, 
{\it Mathematica} will fail to solve the differential equation
produced). 

{\small
\begin{verbatim}
In[50]:= conv[4] = Convert[Sum[ChebyshevT[k,x] z^k, {k, 0, Infinity}], z]

ps-info: PowerSeries  version 1.02, Mar 07, 1994
ps-info: 2 step(s) for RE.
         a[-2 + k] - 2*x*a[-1 + k] + a[k] == 0
ps-info: DE: 
         2 F[z] + 4*(-x + z) F'[z] + (1 - 2*x*z + z^2)\
 
>    F''[z] == 0
ps-info: Trying to solve DE ...
ps-info: DSolve computes 
         (z*C[1] - C[2] + (z*C[2])/C[1])/
 
>    (-1 + 2*x*z - z^2)
ps-info: expression rearranged:
         (z*C[1] + C[2])/(-1 + 2*x*z - z^2)
ps-info: Calculation of initial values...
ps-info: -C[2] == 1
ps-info: (-2*C[1] - 4*x*C[2])/2 == x

                  1                 x z
Out[50]= -(---------------) + ---------------
                         2                  2
           -1 + 2 x z - z     -1 + 2 x z - z

In[51]:= conv[5] = Convert[Sum[LaguerreL[k,a,x] z^k, {k, 0, Infinity}], z]

ps-info: PowerSeries  version 1.02, Mar 07, 1994
ps-info: 2 step(s) for RE.
         (-1 + a + k)*a[-2 + k] + 
 
>     (1 - a - 2*k + x)*a[-1 + k] + k*a[k] == 0
ps-info: DE: 
         (-1 - a + x + z + a*z) F[z] + ((-1 + z)^2)\
 
>    F'[z] == 0
ps-info: Trying to solve DE ...
ps-info: DSolve computes 
         E^(x/(-1 + z))*(-1 + z)^(-1 - a)*C[1]
ps-info: expression rearranged:
         E^(x/(-1 + z))*(1 - z)^(-1 - a)*C[1]
ps-info: Calculation of initial values...
ps-info: C[1]/E^x == 1

          (x z)/(-1 + z)        -1 - a
Out[51]= E               (1 - z)
\end{verbatim}
}\noindent

\section*{Calculation of Recurrence Equations}
\label{sec:Examples on the calculation of recurrence equations}

As described above,
the function {\tt Convert} works 
by generating a 
recurrence equation for the given 
coefficients, coverting it to a differential equation, 
and then solving. The procedure for generating the 
recurrence equation 
is implemented in the function 
{\tt FindRecursion}. 
We give some examples:

{\small
\begin{verbatim}
In[52]:= re[7] = FindRecursion[(1 + (-1)^n)/n, n]

Out[52]= (2 - n) a[-2 + n] + n a[n] == 0

In[53]:= re[8] = FindRecursion[n + (-1)^n, n]

Out[53]= (1 - 2 n) a[-2 + n] - 2 a[-1 + n] +

>     (-3 + 2 n) a[n] == 0

In[54]:= re[9] = FindRecursion[(n + (-1)^n)/n^2, n]

                        2      3
Out[54]= (4 - 12 n + 9 n  - 2 n ) a[-2 + n] +

                     2                   2      3
>     (-2 + 4 n - 2 n ) a[-1 + n] + (-3 n  + 2 n ) a[n]

>      == 0

In[55]:= re[10] = FindRecursion[1/(2n+1)!, n]

Out[55]= a[-1 + n] - 2 n (1 + 2 n) a[n] == 0

In[56]:= ps[8][[1]]

           k       1       k        -2      k   k
         (2  (------------)  - (-----------) ) x
              -1 + Sqrt[5]      1 + Sqrt[5]
Out[56]= ----------------------------------------
                         Sqrt[5]


In[57]:= re[11] = FindRecursion[%, k]

            2
Out[57]= -(x  a[-2 + k]) - x a[-1 + k] + a[k] == 0

In[58]:= re[12] = FindRecursion[E^(-x) LaguerreL[n, alpha, 2x], n]

Out[58]= (-1 + alpha + n) a[-2 + n] +

>     (1 - alpha - 2 n + 2 x) a[-1 + n] + n a[n] == 0

In[59]:= re[13] = FindRecursion[n*LaguerreL[n, alpha, 2x], n]

                                       2
Out[59]= (1 - alpha - 2 n + alpha n + n ) a[-2 + n] +

                                         2
>     (-2 + 2 alpha + 5 n - alpha n - 2 n  - 4 x +

                                        2
>        2 n x) a[-1 + n] + (2 - 3 n + n ) a[n] == 0

In[60]:= re[14] = FindRecursion[LaguerreL[n, alpha, 2x]/n, n]

                                         2
Out[60]= (2 - 2 alpha - 3 n + alpha n + n )

>      a[-2 + n] + (-1 + alpha + 3 n - alpha n -

            2                             2
>        2 n  - 2 x + 2 n x) a[-1 + n] + n  a[n] == 0

In[61]:= re[15] = FindRecursion[LaguerreL[n, x]^2, n]

                        2      3                  2
Out[61]= (4 - 12 n + 9 n  - 2 n  + 4 x - 4 n x + n  x)

                                    2      3
>      a[-3 + n] + (-6 + 22 n - 21 n  + 6 n  - 14 x +

                      2        2        2    3
>        26 n x - 11 n  x - 7 x  + 6 n x  - x )

                                   2      3
>      a[-2 + n] + (2 - 10 n + 15 n  - 6 n  + 6 x -

                      2        2        2    3
>        18 n x + 11 n  x + 5 x  - 6 n x  + x )

                        2      3    2
>      a[-1 + n] + (-3 n  + 2 n  - n  x) a[n] == 0
\end{verbatim}
}\noindent
% COMMENT In[53] braucht ebenfalls laenger...

\vspace*{5mm}
{\Large \bf Electronic Supplement}\abs
The electronic supplement contains the package {\tt PowerSeries.m}.

\end{document}